\documentclass[graybox]{svmult}

\usepackage{mathptmx}       % selects Times Roman as basic font
\usepackage{helvet}         % selects Helvetica as sans-serif font
\usepackage{courier}        % selects Courier as typewriter font
\usepackage{type1cm}        % activate if the above 3 fonts are
\usepackage{makeidx}         % allows index generation
\usepackage{graphicx}        % standard LaTeX graphics tool
\usepackage{multicol}        % used for the two-column index
\usepackage[bottom]{footmisc}% places footnotes at page bottom
\usepackage{amsmath,amssymb,amsmath}
\usepackage{bm}
\usepackage{color}

\makeindex             % used for the subject index
\newcommand{\Rset}{\mathbb R}
\newcommand{\vect}[1]{\boldsymbol{#1}}

\newcommand{\0}{\vect 0}
\newcommand{\bl}{\vect \lambda}
\newcommand{\bx}{\vect \xi}
\newcommand{\f}{\vect f}
\newcommand{\x}{\vect x}
\newcommand{\bv}{\vect v}
\newcommand{\bu}{\vect u}
\newcommand{\z}{\vect z}
\newcommand{\za}{\vect a}
\newcommand{\wb}{\vect b}

\begin{document}

\title*{On the numerical stability of discretised Optimal Control Problems}
\author{Ashutosh Bijalwan and José J Muñoz}
\institute{Ashutosh Bijalwan \at Universitat Politècnica de Catalunya (UPC), Jordi Girona 31, 08034 Barcelona, Spain, \email{abijalwan@cimne.upc.edu}
\and José J Muñoz \at Universitat Politècnica de Catalunya (UPC), Jordi Girona 31, 08034 Barcelona, Spain, 
\email{j.munoz@upc.edu}}
\maketitle

\abstract{Optimal Control Problems consist on the optimisation of an objective functional subjected to a set of Ordinary Differential Equations. In this work, we consider the effects on the stability of the numerical solution when this optimisation is discretised in time. In particular, we analyse a OCP with a quadratic functional and linear ODE, discretised with \textit{Mid-point} and \textit{implicit Euler}. We show that the numerical stability and the presence of numerical oscillations depends not only on the time-step size, but also on the parameters of the objective functional, which measures the amount of control input. Finally, we also show with an illustrative example that these results also carry over non-linear optimal control problems.}

%_______________________________________________________________________________________________________________________________
%________________________________________________________     New Section      _______________________________________________________
%_______________________________________________________________________________________________________________________________

\section{Introduction}
\label{sec:1}
Optimal control is a class of mathematical optimisation problems where the desired system state is determined by minimising/maximising a cost functional subjected to path constraints, written as ordinary differential equations (ODEs) and initial conditions. In a real-world applications, close-form solutions of these problems are difficult to obtain and often they are solved numerically with non-linear programming techniques \cite{bryson75, conway12}. Based on the sequence of optimisation and time discretisation, solution procedures for Optimal Control Problems (OCPs) can be classified as indirect and direct approaches. The indirect method first derives the necessary optimality conditions and forms the Differential-Algebraic Equations (DAEs) with the two-point boundary conditions, also known as the Two-Point Boundary Value Problem (TPBVP) or Hamiltonian Boundary-Value Problem (HBVP), which are then discretisated in time, yielding a discrete system \cite{miller15, sharp2021}. Conversely, the direct method first introduces the time discretisation of the continuous system and then the necessary optimality conditions, derived from the resulting discrete system \cite{Bet10}. Both approaches have their unique properties and limitations, which have been described elsewhere, e.g.  \cite{conway12}.

The stability of numerical integrators for Initial Value Problem (IVP) has been extensively studied. For instance, \emph{explicit Euler (eE)} integrator is only conditionally stable whereas symplectic integrator viz. \emph{implicit Euler (iE)} or \emph{Mid-Point (MP)} are unconditionally stable schemes, which preserve first-integrals of motion \cite{gonzalez99, HWL06}. Due to these promising properties of symplectic integrators, HBVP are frequently discretised with MP or iE schemes \cite{betsch17, flasskamp19, koch13}. However, due to the presence of two-point boundary conditions, implicit (symplectic) integration, which is unconditionally stable along the time-marching direction (forward direction), becomes  explicit (conditionally stable) in the reverse direction, which the natural direction of the adjoint ODEs in OCP. This fact motivates the study in this work, where we show that originally stable schemes for ODEs, may become unstable in discretised OCP. Furthermore, we show that the stability and the presence of oscillation in the numerical solution depend on the magnitude of the control input. For simplicity and clarity on the exposition of our ideas, we have chosen to restrict our discussion to the numerical stability of MP and iE schemes applied to OCPs. 

The chapter is structured as follows. In Section~\ref{sec:1}, we begin with a brief review of the available numerical methods for the OCP. In Section~\ref{sec:2}, we formally introduce the general framework of continuous OCP and describe the associated symplectic structure. Additionally, we present the time-discrete form of the HBVP. In Section~\ref{sec:3}, we employ the common MP and iE schemes on an illustrative  linear OCP, and provide a detailed analysis of the source of numerical oscillations and associated stability criteria. Finally, in Section~\ref{sec:4}, we generalise these ideas to prevent numerical oscillations in the non-linear OCPs, and conclusions are drawn in  Section~\ref{sec:5}.  

%_______________________________________________________________________________________________________________________________
%________________________________________________________     New Section      _______________________________________________________
%_______________________________________________________________________________________________________________________________

\section{Optimal Control Problem: Indirect Method}
\label{sec:2}
Let us consider a continuous optimal control problem which seeks control $\bu(t)$ and state $\x(t)$ trajectories  that minimise an objective functional $\mathcal J(\x, \bu)$ subjected to first-order IVP and equality constraints \cite{Bet10}, i.e. are solution of the following optimisation problem
\begin{align}\label{e:1} 
\max_{\x(t),\text{ }\bu(t) } & \mathcal{J}(\x, \bu) \\ \nonumber
\text{s.t.},&\ \f(\x,\bu)-\dot{\x}=\0 && \text{(State ODE)}\\ \nonumber
&\ \x(0)-\x_0 =\0 &&\text{(Initial conditions)}
\end{align}
where for clarity, we suppress time-dependence. A common choice for the functional $\mathcal J(\x, \bu)$ in trajectory optimisation problems is to provide a measure of deviation of the state variable $\x(t)$ from a target state $\x_t$, and also add an associated input/control cost to achieve the desired state. Mathematically, the following form is usually employed,
\begin{align}\label{e:2} 
\mathcal{J}(\x, \bu):=\int^T_0 \left(\frac{1}{2}(\x-\x_t)^\mathsf{T} \textbf{R} (\x-\x_t)+\frac{\alpha}{2}\bu^\mathsf{T} \textbf{Q} \bu \right) dt,
\end{align}
where $\textbf{Q}$ and $\textbf{R}$ are the input and output matrices. Parameter $\alpha >  0$ regulates the amount of input control $\bu$ and $T$ is the fixed final time.

In order to deduce the optimality conditions of the the constrained optimisation problem in Eq. \eqref{e:1}, we introduce the Lagrangian multipliers $\bl(t)$ and $\bx$ and define the Lagrangian functional associated with  problem \eqref{e:1} as
\begin{align}\label{e:3}
\mathcal{L}(\x, \bu; \bl, \bx)=\mathcal{J}(\x, \bu)+\int^T_0 \bl^\mathsf{T} \left(\f(\x,\bu)-\dot{\x}\right) dt +\bx^\mathsf{T} (\x(0)-\x_0). 
\end{align}

After integrating by parts and rearranging the integrands, the Lagrangian becomes \cite{bryson75}
\begin{align*}
\mathcal{L}(\x, \bu; \bl, \bx)=\int^T_0 \left(\mathcal{H}+\dot{\bl}^\mathsf{T}\x \right) dt-\bl(T)^\mathsf{T}\x(T) +\bl(0)^\mathsf{T}\x(0)+\bx^\mathsf{T}(\x(0)-\x_0)  
\end{align*} 
where $\mathcal{H}(\x, \bu; \bl):=\frac{1}{2}(\x-\x_t)^\mathsf{T} \textbf{R} (\x-\x_t)+\frac{\alpha}{2}\bu^\mathsf{T} \textbf{Q} \bu + \bl^\mathsf{T} \f(\x,\bu)$ is the so-called Control Hamiltonian, which is conserved along optimal trajectories, provided $\frac{\partial \mathcal{H}}{\partial t}=0$. First-order optimality conditions can be then derived from the derivatives of $\mathcal L$ with respect to $\x, \bl, \bu$ and $\bx$, which yields the following system of differential-algebraic equations (DAEs) with the two-point boundary conditions \cite{betsch17, flasskamp19},
\begin{subequations}\label{e:5}
\begin{align}
\dot{\bl}=& -\nabla_{\x} \mathcal{H} \\
\dot{\x}=& \nabla_{\bl} \mathcal{H}\\
\0=& \nabla_{\bu} \mathcal{H} \\
\x(0)=& \x_o \text{,  } \bl(T)=\0\label{e:5ic}
\end{align}
\end{subequations}

The first three equations above are commonly named adjoint, state, and control-algebraic equations, respectively,  and the system is known as the set of Euler-Lagrange equations (E-L), which together with the end conditions in \eqref{e:5ic} constitutes the so-called Hamiltonian boundary-value problem (HBVP),  with inherent symplectic structure \cite{bryson75}. 

We resort now to the numerical time integration of the HBVP by introducing a time discretisation scheme of the continuous system in Eq. \eqref{e:5} with a uniform time-step size $\Delta t >0$ ($n=1,2,\dots,N$), which we write in the following general form:
\begin{subequations}\label{e:6}
\begin{align}
\frac{\bl_{n-1}-\bl_n}{\Delta t} &=\textbf{R}\left(\x_{n-\tau}-\x_t\right)+ \nabla_{\x} \f (\x_{n-\tau},\bu_{n-\tau}) ^\mathsf{T} \bl_{n-\tau} \\
\frac{\x_n-\x_{n-1}}{\Delta t} &=\f (\x_{n-\tau},\bu_{n-\tau}) \\
\0 &=\alpha \textbf{Q}\bu_{n-\tau}+\nabla_{\bu} \f (\x_{n-\tau},\bu_{n-\tau}) ^\mathsf{T} \bl_{n-\tau} \\
\x_0 &=\x_o \text{, } \bl_N=\0
\end{align}
\end{subequations}
with $\x_{n-\tau}:=\tau \x_{n-1}+(1-\tau) \x_n$, and similarly for $\bu_{n-\tau}$ and $\bl_{n-\tau}$. The values $\tau=0$ and $\tau=\frac{1}{2}$ correspond to the well known implicit Euler (iE) and Mid-Point (MP) schemes, respectively.

%_______________________________________________________________________________________________________________________________
%________________________________________________________     New Section      _______________________________________________________
%_______________________________________________________________________________________________________________________________

\section{Motivating Example: Linear OCP}
\label{sec:3}
Let us consider a propelled body (e.g. jellyfish) with mass $m$ moving along $+y$ direction starting with an initial velocity $v_o$. There is a gravitational force with acceleration $a$ acting along the $-y$ direction,  in addition to a drag force with the form $f_d:= b v $, with $b>0$ a frictional coefficient. Furthermore, we assume that the jellyfish controls the propulsion velocity by regulating the fluid ejection along $-y$ direction resulting in a thrust force $u$ along $+y$ direction. Linear momentum balance along $+y$ direction results in the following state IVP,
\begin{eqnarray}
f(v,u)-\dot{v} & =0   \nonumber \\
v(0)-v_o &=0
\label{e:8}
\end{eqnarray}
where $f(v,u):=-\frac{b}{m}v+\frac{u}{m}-a$. We are interested in obtaining a control function $u(t)$, which drives the jellyfish from initial velocity $v_o$ to a target velocity $v_t$ in a given time $T$. Essentially, we want to minimise the following cost functional 
\begin{align}\label{e:9} 
\mathcal{J}(v, u):=\int^T_0 \left(\frac{1}{2}(v-v_t)^2 +\frac{\alpha}{2} u^2 \right) dt
\end{align}
subject to the IVP in Eq. \eqref{e:8}. A closed-form solution of the posed problem is given by
\begin{subequations}\label{e:10}
\begin{align}
v(t)  &= C_1 \left(\frac{b}{m}-\gamma \right) e^{\gamma t}+C_2  \left(\frac{b}{m}+\gamma \right) e^{-\gamma t}+v_p   \label{e:2c58a} \\
\lambda(t)  &= C_1 e^{\gamma t}+C_2 e^{-\gamma t}+\lambda_p  \label{e:2c58b} \\
u(t) &= \frac{\lambda(t)}{\alpha m}
\end{align}
\end{subequations}
where 
\begin{align*}
&\gamma=\sqrt{\frac{b^2}{m^2}+\frac{1}{\alpha m^2}}, \ v_p=\frac{v_t-\alpha b m a }{\alpha m^2 \gamma^2}, \ \lambda_p=-\frac{b v_t+m a }{ m \gamma^2}, \\ &
C_1=\frac{m \left(v_o-v_p\right) e^{-\gamma T}+ \left(b+m \gamma \right) \lambda_p}{ \left(b-m \gamma \right) e^{-\gamma T}  -  \left(b+m \gamma \right) e^{\gamma T}}, \ C_2=\frac{m \left(v_o-v_p\right)- \left(b-m \gamma \right) C_1}{ b+m \gamma}
\end{align*}

Irrespective of the value of $\alpha>0$, we have that $\gamma>0$ and one would expect a non-oscillatory state, adjoint and control functions. We will show in the next sections that the numerical solution of the OCP with the MP and iE time discretisation schemes is not necessarily non-oscillatory. 

%________________________________________________________    Subsection _______________________________________________________

\subsection{Mid-Point (MP) scheme}
\label{subsec:3}
By using the value $\tau=1/2$ in Eq. \eqref{e:6}, a mid-point (MP) time discretisation of the ODE system is obtained,
\begin{subequations}\label{e:12}
\begin{align}
\frac{\lambda_n-\lambda_{n-1}}{\Delta t}-\frac{b}{m}\lambda_{n-\frac{1}{2}}+v_{n-\frac{1}{2}}-v_t&=0, \label{e:12La} \\
\frac{v_n-v_{n-1}}{\Delta t}+\frac{b}{m} v_{n-\frac{1}{2}}-\frac{1}{m} u_{n-\frac{1}{2}}+a&=0, \label{e:12Lb} \\ 
u_{n-\frac{1}{2}}+\frac{1}{\alpha m} \lambda_{n-\frac{1}{2}}&=0. \label{e:12Lc}
\end{align}
\end{subequations}

After substituting the boundary conditions $v_0=v_o$ and $\lambda_N=0$, the resulting linear system of equations can be solved with conventional linear solvers.

In order to study the stability of the the MP scheme, it will be helpful to express the system with the independent variables ($v$ and $\lambda$). Notice that the control equation in Eq. \eqref{e:12Lc} is linear, so that replacing Eq. \eqref{e:12Lc} into Eq. \eqref{e:12Lb}, with the definitions
\begin{align}
p:=\frac{b}{m}+\frac{2}{\Delta t} ; \   q:=\frac{b}{m}-\frac{2}{\Delta t} ; \  s:=\frac{1}{\alpha m^2}
\end{align}
results in the following reduced discrete system 
\begin{subequations}\label{e:13}
\begin{align}
& q \lambda_n+p \lambda_{n-1}-v_{n-1} -v_n +2 v_t=0 \\
& p v_n+q v_{n-1} + s \lambda_{n-1}+ s \lambda_n+2a=0
\end{align}
\end{subequations}

By defining the vector $\z_n:=\{v_n,\lambda_n\}^T$, the MP scheme can be expressed as
\begin{align}\label{e:14}
\z_n = \textbf{A}  \z_{n-1} + \za
\end{align}
where 
\begin{align*}
\mathbf A:=\frac{-1}{s+p q} \begin{bmatrix}
s+q^2   & -s (p-q)    \\
-(p-q) &  s+p^2     \\
\end{bmatrix}; \ \za:= \frac{-2}{s+p q} \left\{\begin{array}{c}
q a  -s v_t      \\
a+ p v_t   
\end{array}\right\}.     
\end{align*}

The eigenvalues of matrix $\textbf{A}$  are real and distinct, and are given by 
\begin{align}\label{e:15}
e_{1,2} = \left \{ 
\frac{2+  \gamma \Delta t }{2-  \gamma \Delta t}, \ 
\frac{2-  \gamma \Delta t }{2+ \gamma \Delta t}  
\right \}.
\end{align}
and the spectral radius of $\textbf{A}$ is 
\begin{align}\label{e:16}
\rho (\mathbf A)=\frac{2+ \gamma \Delta t }{\left|  2- \gamma \Delta t  \right|}. 
\end{align}

Consequently, the time discretisation scheme is stable if $\rho(\textbf{A}) \leq 1$. It seems that there exist no admissible pair ($\gamma$, $ \Delta t) \in \Rset^+$ for which $\rho(\textbf{A}) \leq 1$. However, if we impose restrictions on $\gamma \Delta t $ such that $\log_{10} (\gamma \Delta t/2) \in (-\infty, -\epsilon) \cup (\epsilon, \infty) $ with $ \epsilon \approx 2$, numerical evidence shows that the solution remains bounded.

In order to analyse the presence of numerical oscillations, we examine the value of $\alpha$ for which the eigenvalues change their sign. The eigenvalues in Eq. \eqref{e:15}
switch sign from negative to positive if $\gamma^2-\frac{4}{ \Delta t^2 }  < 0$ and for the threshold value $\gamma=\frac{2}{\Delta t} $, system blows out. In summary,
\begin{align}\label{e:17}
e_i =\begin{cases}
<0 , & \alpha < \alpha_{th,MP} \ \text{ (Oscillatory response)} \\
> 0, & \alpha > \alpha_{th,MP} \ \text{ (No oscillations)}
\end{cases}
\end{align}
with
\begin{align}\label{e:aMP}
\alpha_{th,MP}:=\frac{\Delta t^2}{4 m^2-b^2 \Delta t^2}.
\end{align}

%________________________________________________________    Subsection _______________________________________________________

\subsection{Implicit Euler (iE) scheme}\label{subsec:3b}
By using the value $\tau=0$ in Eq. \eqref{e:6}, Implicit Euler (iE) time discretisation results in the following discrete system 
\begin{subequations}\label{e:19}
\begin{align}
\frac{\lambda_n-\lambda_{n-1}}{\Delta t}-\frac{b}{m}\lambda_{n}+v_{n}-v_t&=0, \label{e:19La} \\
\frac{v_n-v_{n-1}}{\Delta t}+\frac{b}{m} v_{n}-\frac{1}{m} u_{n}+a&=0, \label{e:19Lb} \\
u_{n}+\frac{1}{\alpha m} \lambda_{n}&=0. \label{e:19Lc}
\end{align}
\end{subequations}

After substituting the boundary conditions $v_0=v_o$ and $\lambda_N=0$, the solution of the resulting linear system can be obtained. 

In order to study the stability of the iE scheme, we express our system in terms of the independent variables $v$ and $\lambda$. Substituting Eq. \eqref{e:19Lc} into Eq. \eqref{e:19Lb}, and with the new definitions
\begin{align}
p^*:=\frac{b}{m}+\frac{1}{\Delta t} ; \  q^*:=\frac{b}{m}-\frac{1}{\Delta t} ; \ r:=\frac{1}{\Delta t}; \ s:=\frac{1}{\alpha m^2}
\end{align}
the system in \eqref{e:19} is equivalent to
\begin{align}\label{e:20}
& q^* \lambda_n+r \lambda_{n-1} -v_n + v_t=0 \\
& p^* v_n-r v_{n-1} + s \lambda_n+a=0.
\end{align}
or in terms of the vector $\z_n:=\{v_n,\lambda_n\}^T$, it canbe expressed as 
\begin{align}\label{e:21}
\z_n = \textbf{B}  \z_{n-1} + \wb 
\end{align}
with
\begin{align*}
\mathbf B:=\frac{-r}{s+p^* q^*} \begin{bmatrix}
-q^*   & -s     \\
-1 &  p^*     \\
\end{bmatrix}, \ \wb:= \frac{-1}{s+p^* q^*} \left\{\begin{array}{c}
q^* a  -s v_t      \\
a+ p^* v_t       
\end{array}\right\}.
\end{align*}

The eigenvalues of matrix $\mathbf B$ and its spectral radius $\rho(\mathbf B)$ are given by
\begin{align}\label{e:22}
e_{1,2}=\left \{ \frac{1}{1 + \gamma \Delta t}   ,  \frac{1}{1 - \gamma \Delta t} \right \}, \  \rho (\mathbf B)=\frac{1}{\left|  1- \gamma \Delta t \right| }
\end{align}

The time discretisation scheme will be thus stable if $\left| \gamma \Delta t-1 \right|  \geq 1 $. Since $ \gamma \Delta t \in \Rset^+$, we are left with the restriction $\gamma \geq \frac{2}{\Delta t}$.

\begin{figure}[!htb]
    \centering   
    \includegraphics[width=0.495\textwidth]{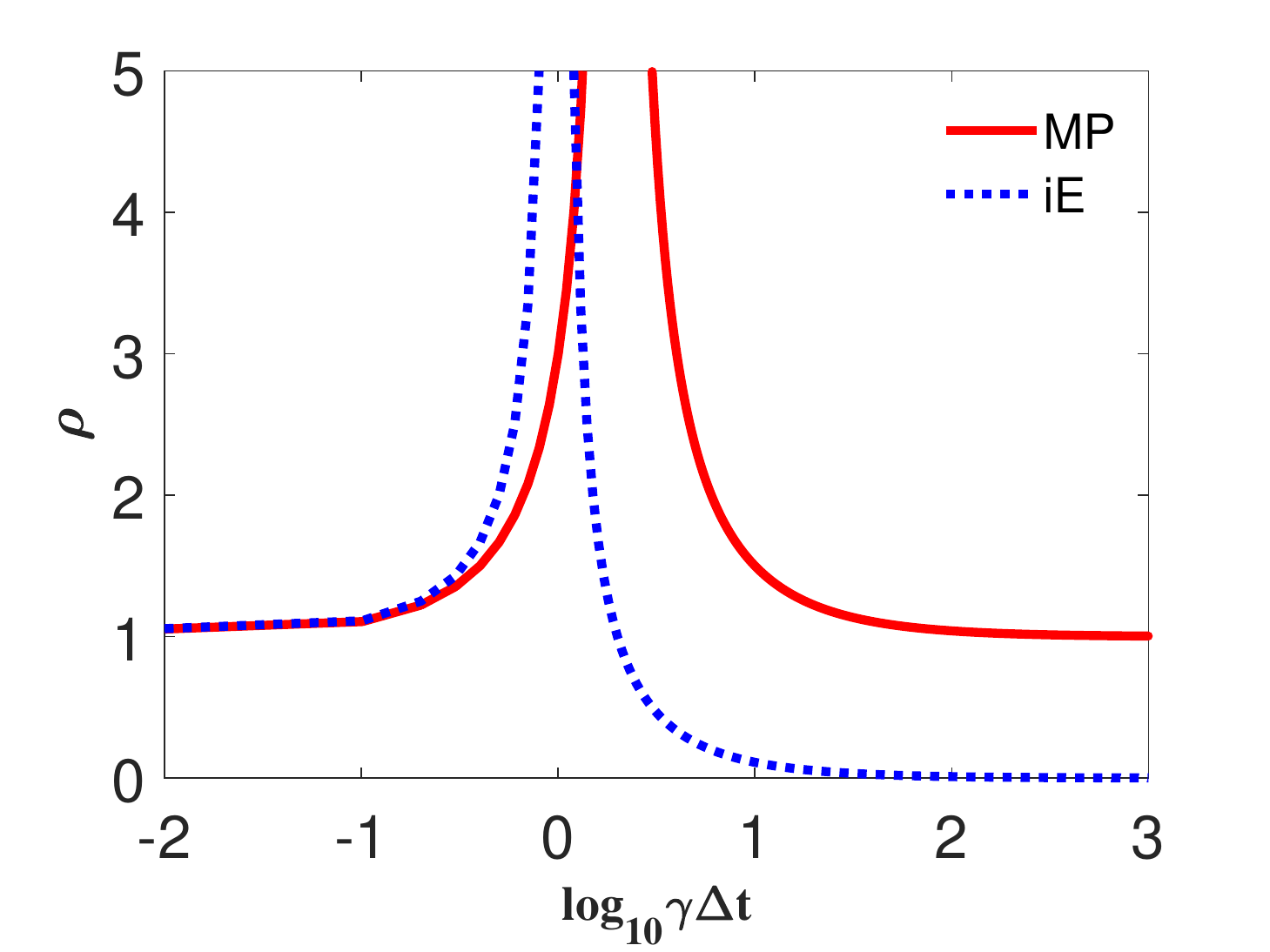} 
    \includegraphics[width=0.495\textwidth]{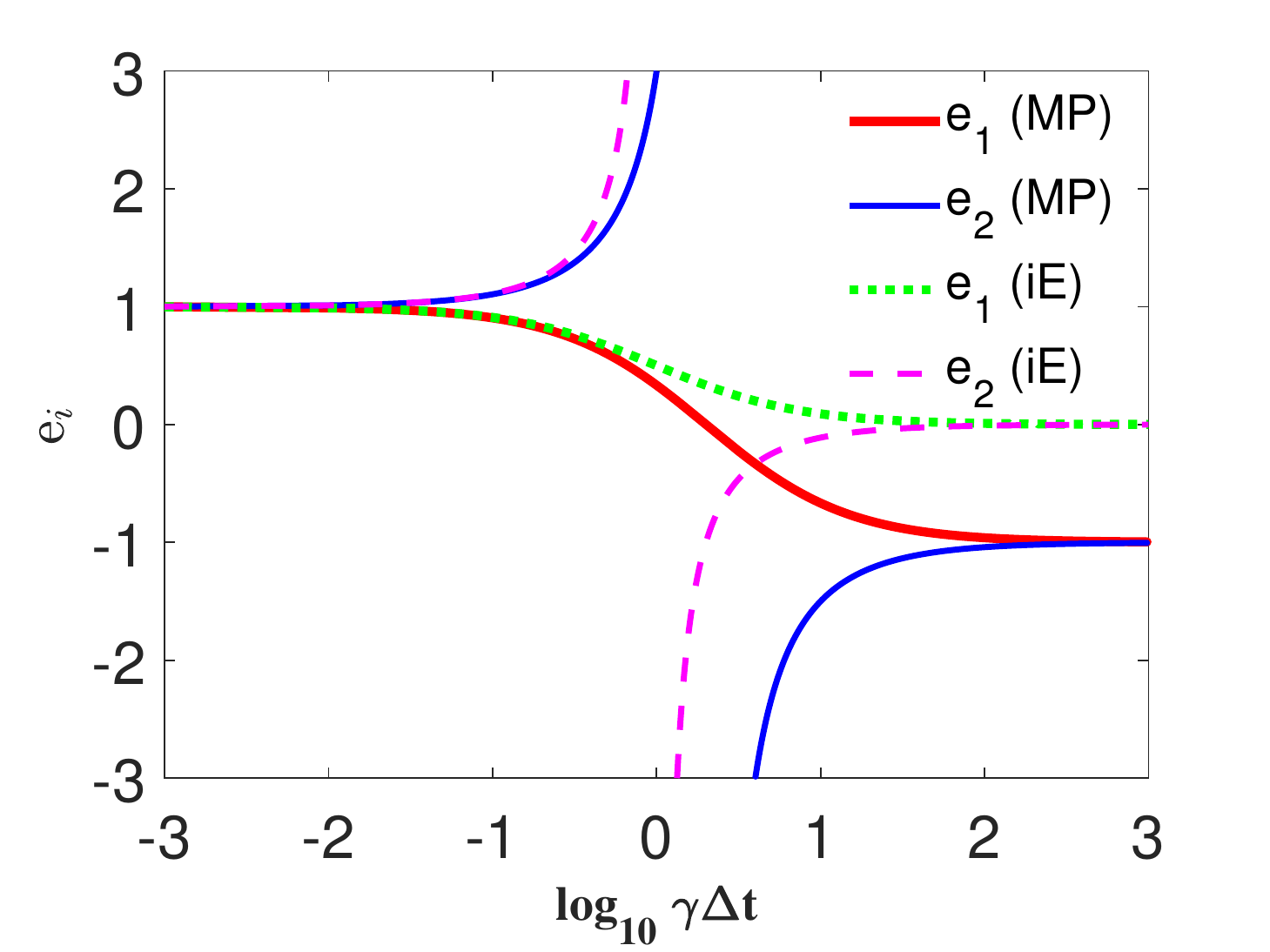} 
    \caption{Mid-point (MP) and implicit Euler (iE) scheme: (a) spectral radius, and (b) eigenvalues.}
   \label{fig:1}
\end{figure}

\begin{figure}[!htb]
    \centering   
    \includegraphics[width=0.495\textwidth]{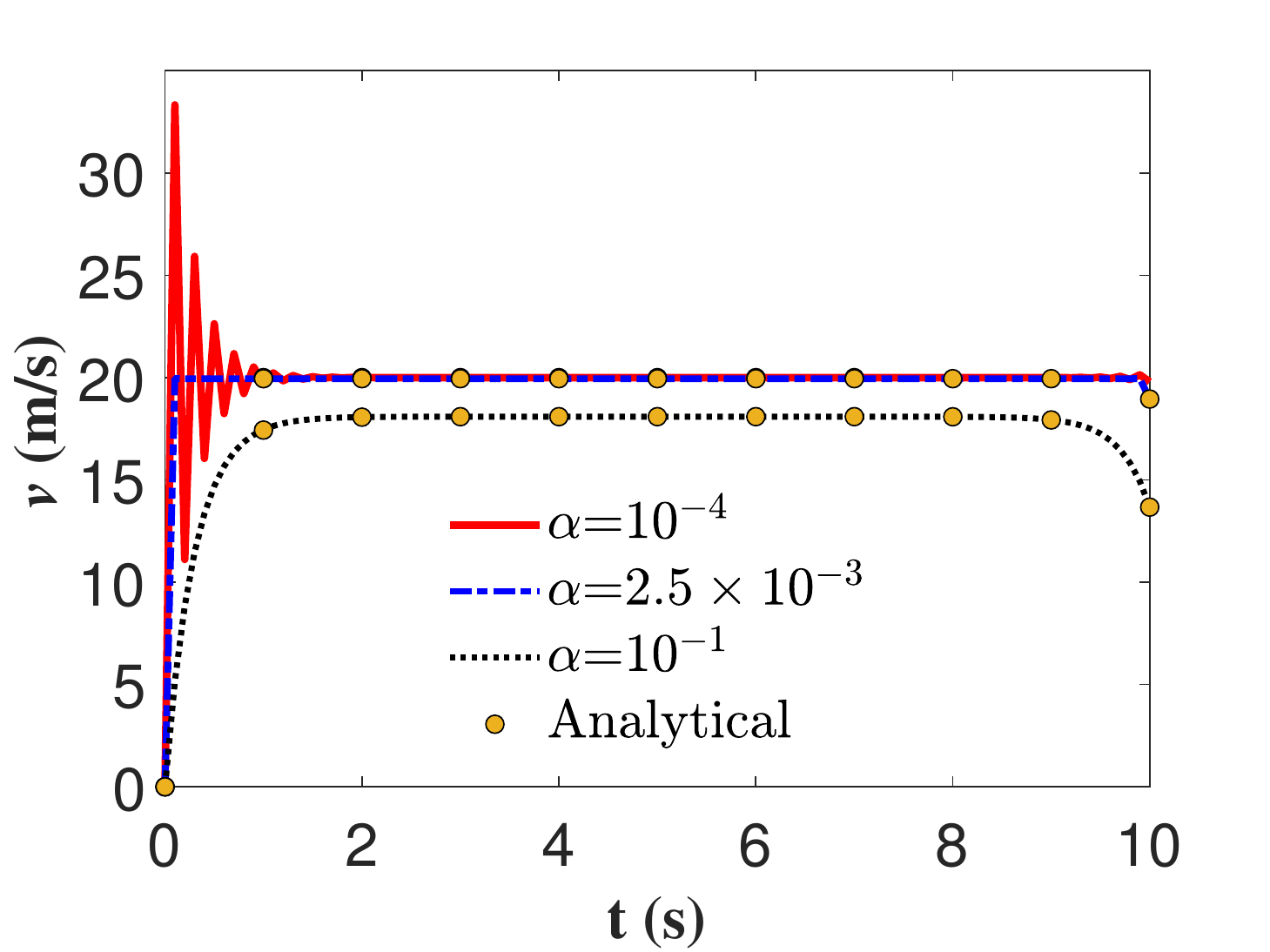} 
    \includegraphics[width=0.495\textwidth]{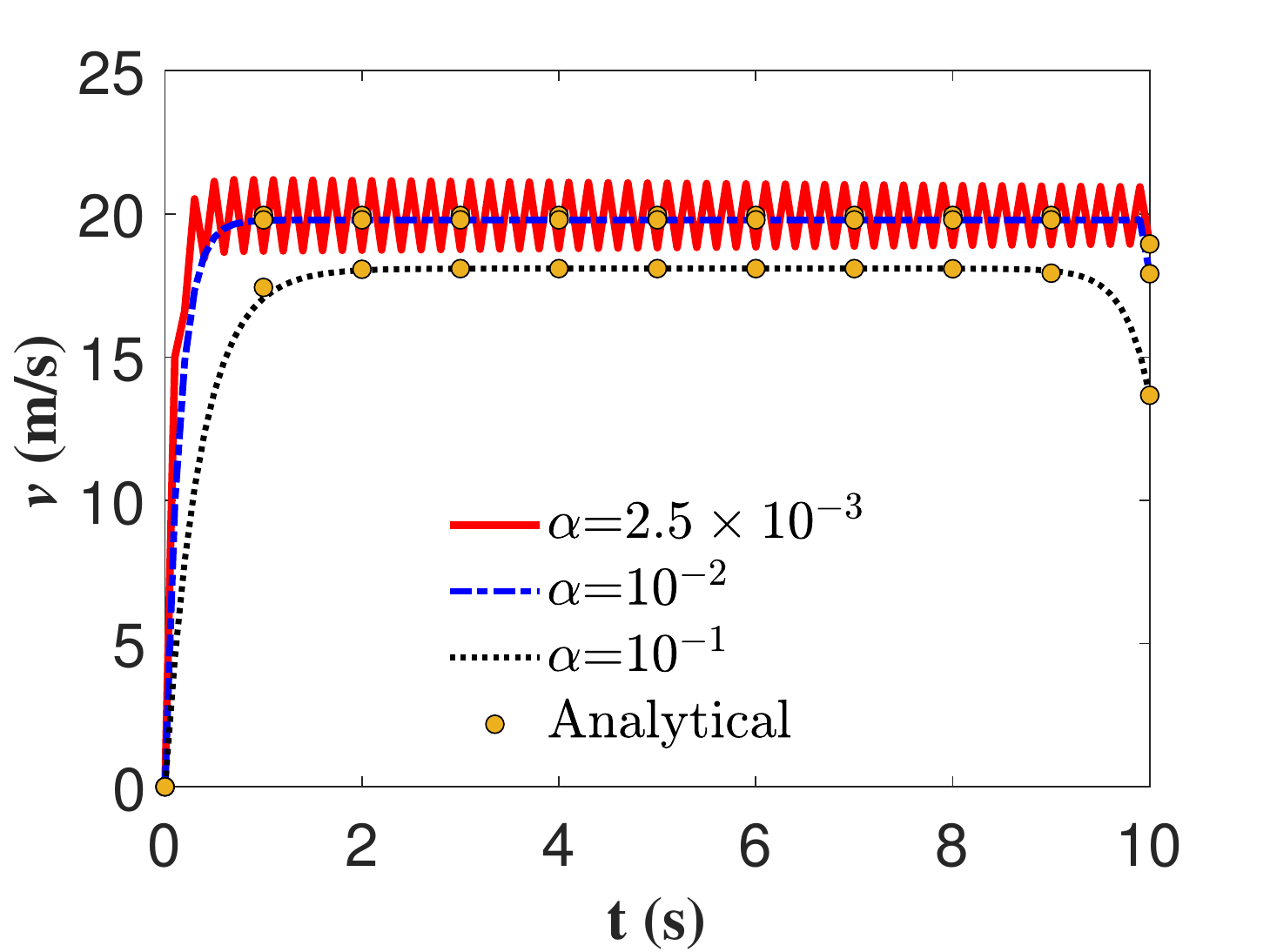}
    \caption{State optimal velocity: (a) MP scheme and (b) iE scheme (dots: Analytical solution).} 
   \label{fig:2}
\end{figure}

\begin{figure}[!htb]
    \centering   
    \includegraphics[width=0.495\textwidth]{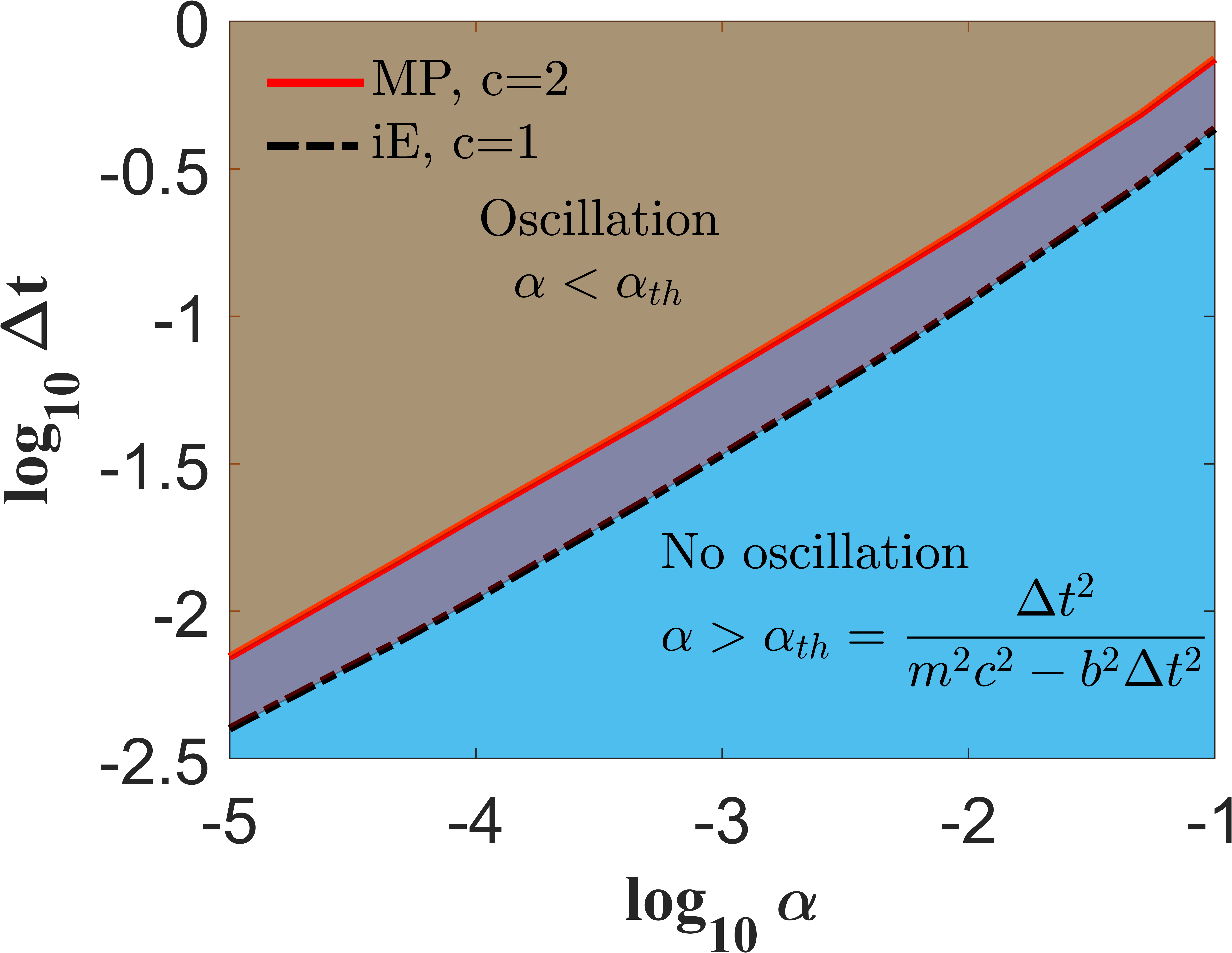} 
    \caption{Phase diagram $(\alpha, \Delta t)$ for MP and iE schemes.}
   \label{fig:3}
\end{figure}

A similar analysis to the MP scheme inidcates that the change of sign in  $e_2$ in the iE scheme occurs when $1- \gamma \Delta t$ changes sign, and for the threshold value of $\gamma=\frac{1}{\Delta t} $ the system blows out. Summarising,
\begin{align}\label{e:24}
\frac{1}{1- \gamma \Delta t} =\begin{cases}
< 0, & \alpha < \alpha_{th,iE} \text{ (Oscillatory response)}\\
> 0, & \alpha > \alpha_{th,iE} \text{ (No oscillations)}
\end{cases}
\end{align}
with
\begin{align}\label{e:aiE}
\alpha_{th,iE}:=\frac{\Delta t^2}{ m^2-b^2 \Delta t^2}.
\end{align}

We have numerically verified the theoretical results and thresholds in \eqref{e:aMP} and \eqref{e:aiE} by using the values $(m,b,a,v_o,v_t,T,\Delta t)=(1,1,1,0,20,10,0.1)$, which imply the values $\alpha_{th,MP}=2.5E-3$ and $\alpha_{th,iE}=1.01E-2$. Reducing $\alpha$ below these values, numerical oscillations are obtained, in agreement with conditions in  \eqref{e:17} and \eqref{e:24}, as shown in Fig. \ref{fig:2}. Theoretical results for the spectral radius and eigenvalues are shown in Fig.~\ref{fig:1}. We point out that there exists a region where both schemes become unstable and that stabile results are obtained if we are sufficiently far from this region.

For verifying the thresholds of  $\alpha_{th}$, we generate the stability envelope from numerical experiments in Fig.~\ref{fig:3}. We found close matches with Eq. \eqref{e:aMP} and Eq. \eqref{e:aiE}. We remark that the phase boundary of the MP scheme lies above the iE scheme and shows numerical oscillations much later than the iE scheme.

%-----------------------------------------------------------------------------------------------------------------------------------------------------------------
%------------------------------------------------------------------ New Section ------------------------------------------------------------------------------
%-----------------------------------------------------------------------------------------------------------------------------------------------------------------

\section{Nonlinear OCP: Inverted Elastic Pendulum}\label{sec:4}

\begin{figure}[!htb]
    \centering
    \includegraphics[width=0.48\textwidth]{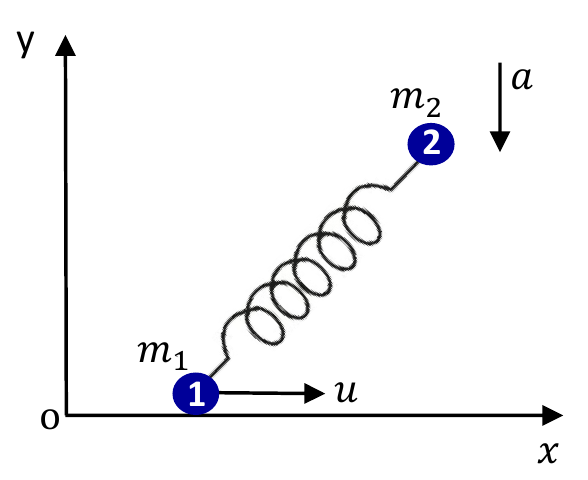} 
    \includegraphics[width=0.51\textwidth]{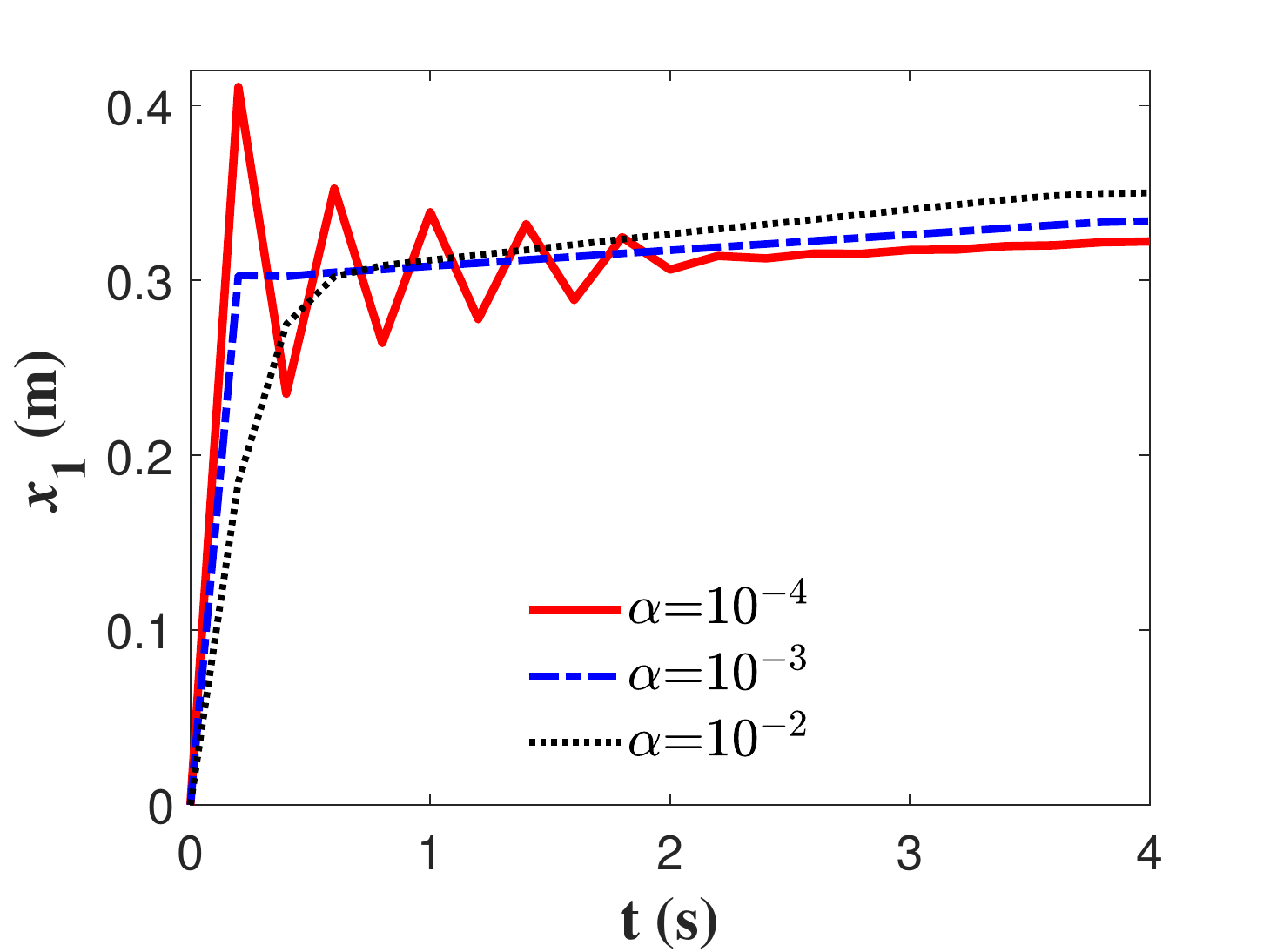}
    \caption{Inverted pendulum: mass $m_1$ horizontal displacement with $\alpha$ (MP Scheme)} 
   \label{fig:4}
\end{figure}

In the second application, we consider a planar elastic inverted pendulum consisting on two point masses $m_1$ and $m_2$ and linked by an elastics spring. The system is subjected to a gravitationsl force field along $-y$ direction with intensity $a$. Position $\x_1$ of $m_1$ is constrained to move along x-axis, while $m_2$ is free to move in x-y plane. The initial positions of $m_1$ and $m_2$ are respectively $(0,0)$ and $(0.3,1)$ (see Fig.~\ref{fig:4}). We assume that the spring potential energy varies quadratically with the length increment from a rest-length $l_o$, i.e. $U(\x):=\frac{k}{2}(l(\x)-l_o)^2$. Linear momentum balance can be expressed as
\begin{align}\label{e:25}
\f(\x, \bu)-\dot{\bv} &=\0 \\ \nonumber
\bv-\dot{\x} &=\0 \\ \nonumber
\x(0) - \x_o &=  \0 \\ \nonumber
\bv(0) - \bv_o &= \0
\end{align}
where $l=\| \x_2-\x_1 \|_2$, $\x=\{\x_1,\ \x_2\}^\mathsf{T}$, $\bv= \{\bv_1,\bv_2\}^\mathsf{T}$,  $\f(\x, \bu):=-\mathbf{M}^{-1} \left( \nabla_{\x} U + \za \right)$,  $\mathbf{M}:=\text{diag}(m_1, m_1, m_2,m_2)$, and $\za:=\{0,a,0,a\}^\mathsf{T}$. We are interested in finding the control velocity $u=v_1$ which stabilises the system in upright configuration. Equivalently, we aim at minimising the cost functional  
\begin{align}\label{e:26} 
\mathcal{J}(\x, \bu):=\int^T_0 \left(\frac{1}{2}(\x-\x_t)^\mathsf{T} \textbf{R} (\x-\x_t)+\frac{\alpha}{2}\bu^\mathsf{T} \textbf{Q} \bu \right) dt
\end{align}
subjected to Eq. \eqref{e:25}, where $\bu=\bv$, $\textbf{Q}=\text{diag}(1,0,0,0)$, $\textbf{R}=\text{diag}(0,0,0,1)$, and $\x_t=\{0,0,0,x_t\}^\mathsf{T}$. 

For the numerical test, we assign system parameters $(m_1, m_2, k, a, x_t, T, \Delta t)=(1, 1, 1, 1, 2, 4, 0.2)$, and we use $\alpha \in (10^{-4} , 10^{-2})$. Next, we discretise our system with MP scheme, i.e. Eq. \eqref{e:6} with $\tau=\frac{1}{2}$). The resulting system of non-linear equations is  numerically solved. Based on the linear analysis, we observed that for the MP scheme with small $\Delta t$, we have that $\alpha_{th} \sim \mathcal{O} (\Delta t^2)$, hereby coined as a \textit{conservative stability criteria}, and reducing $\alpha$ below this value should trigger numerical oscillations. Fig.~\ref{fig:4} shows the horizontal position of mass 1 with three $\alpha$ values. It can be seen that for $\alpha=10^{-2}$ we have smooth trajectory, with $\alpha=10^{-3}$ we incubate small kink, and $\alpha=10^{-4}$ results in an oscillating optimal trajectory.

%_______________________________________________________________________________________________________________________________
%________________________________________________________     New Section      _______________________________________________________
%_______________________________________________________________________________________________________________________________

\section{Summary}
\label{sec:5}
We have shown that the stability criteria in OCP depends not only on the parameters of the ODEs and time-step of the discretisation, but also on the parameters of the objective function. Furthermore, small values of the parameters in the control cost function may also induce numerical oscillations. We have demonstrated the origin of these instabilities and oscillations for a linear problem, and also illustrated numerically how these ideas also carry over problems with non-linear ODEs. Our numerical experiments suggest that as $\alpha$ diminishes, we must reduce the time-step size to circumvent numerical oscillations in optimal trajectories, i.e., $\alpha \sim \mathcal{O} (\Delta t^2)$.
\begin{acknowledgement}
This work is financially supported by the Spanish Ministry of Science and Innovation, under Severo Ochoa program CEX2018-000797-S, and the research project DynAd2, with reference PID2020-116141GB-I00.
\end{acknowledgement}


\begin{thebibliography}{spphys}

\bibitem[1]{Bet10}
  Betts, J.T.: {Practical Methods for Optimal Control and Estimation Using
  Nonlinear Programming}.
  \newblock Society for Industrial and Applied Mathematics (SIAM), 2010.
  
\bibitem[2]{betsch17}
  Betsch, P.; Becker, C.: {Conservation of generalized momentum maps in mechanical optimal
  control problems with symmetry}.
  \newblock International Journal for Numerical Methods in Engineering, Vol.~111(2), 144-175, doi.org/10.1002/nme.54594, 2017.
  
\bibitem[3]{bryson75}
  Bryson, A.E.; Ho, Y.C.: {Applied Optimal Control. Optimization, Estimation and Control}.
  \newblock Taylor {\&} Francis,  1975.

\bibitem[4]{conway12}
 Conway, B.A.: {A survey of methods available for the numerical optimization of continuous dynamic systems}.
  \newblock Journal of Optimization Theory and Applications, Vol.~152(2), 271-306, doi.org/10.1007/s10957-011-9918-z, 2012.

\bibitem[5]{flasskamp19}
  Fla{\ss}kamp, K.; Murphey, T.D.: {Structure-preserving local optimal control of mechanical systems}.
  \newblock Optimal Control Applications and Methods, Vol.~40(2):310-329, doi.org/10.1002/oca.2479, 2019.

\bibitem[6]{gonzalez99}
Gonzalez, O.: {Mechanical systems subject to holonomic constraints: Differential–algebraic formulations and conservative integration}. 
\newblock Physica D: Nonlinear Phenomena, Vol.~132(1-2), 165-174, doi.org/10.1016/S0167-2789(99)00054-8, 1999.

\bibitem[7]{HWL06}
  Hairer, E.; Wanner, G.; Lubich, C.: {Symplectic integration of hamiltonian systems. In Geometric Numerical Integration}.
  \newblock Springer, Berlin, Heidelberg, 2006.

\bibitem[8]{koch13}
 Koch, M.W.; Leyendecker, S.: {Energy momentum consistent force formulation for the optimal control
  of multibody systems}.
  \newblock Multibody System Dynamics, Vol.~29(4):381-401, doi.org/10.1007/s11044-012-9332-9, 2013.

\bibitem[9]{miller15}
 Miller, M.I.; Trouv{\'e}, A.; Younes, L.: {Hamiltonian systems and optimal control in computational anatomy: 100
  years since d'arcy thompson}.
  \newblock  Annual review of biomedical engineering, Vol.~17:447-509, doi.org/10.1146/annurev-bioeng-071114-040601, 2015.

\bibitem[10]{sharp2021}
 Sharp, J.A.; Burrage, K.; Simpson, M.J.: {Implementation and acceleration of optimal control for systems
  biology}.
  \newblock  J. R. Soc. Interface, Vol.~18(181):20210241, doi.org/10.1098/rsif.2021.0241, 2021.


\end{thebibliography}
\end{document}